# ON THE EVENTUAL PERIODICITY OF FRACTIONAL ORDER DISPERSIVE WAVE EQUATIONS USING RBFS AND TRANSFORM

**Hameed Ullah Jan**✉
*Department of Basic Sciences and Islamiat[1]*
*huj@uetpeshawar.edu.pk*

**Marjan Uddin**
*Department of Basic Sciences and Islamiat[1]*

**Irshad Ali Shah**
*Department of Mathematics[2]*

**Salam Ullah Khan**
*Department of Computer Science[2]*

*[1]University of Engineering and Technology Peshawar*
*Khyber Pakhtunkhwa, Pakistan, 25120*

*[2]University of Science and Technology Bannu*
*Khyber Pakhtunkhwa, Pakistan, 28100*

✉Corresponding author

**Abstract**

In this research work, let's focus on the construction of numerical scheme based on radial basis functions finite difference (RBF-FD) method combined with the Laplace transform for the solution of fractional order dispersive wave equations. The numerical scheme is then applied to examine the eventual periodicity of the proposed model subject to the periodic boundary conditions. The implementation of proposed technique for high order fractional and integer type nonlinear partial differential equations (PDEs) is beneficial because this method is local in nature, therefore it yields and resulted in sparse differentiation matrices instead of full and dense matrices. Only small dimensions of linear systems of equations are to be solved for every center in the domain and hence this procedure is more reliable and efficient to solve large scale physical and engineering problems in complex domain.

Laplace transform is utilized for obtaining the equivalent time-independent equation in Laplace space and also valuable to handle time-fractional derivatives in the Caputo sense.

Application of Laplace transform avoids the time steeping procedure which commonly encounters the time instability issues. The solution to the transformed model is then obtained by computing the inversion of Laplace transform with an appropriate contour in a complex space, which is approximated by trapezoidal rule with high accuracy. Also since the Laplace transform operator is linear, it cannot be used to transform non-linear terms therefore let's use a linearization approach and an appropriate iterative scheme. The proposed approach is tasted for some nonlinear fractional order KdV and Burgers equations. The capacity, high order accuracy and efficiency of our approach are demonstrated using examples and results.

## 1. Introduction

Significant advancement to both the theory and application of fractional calculus has been made in the last century. The most essential advantage of fractional calculus is its nonlocal nature and effectiveness in modeling anomalous diffusion that happens in complex system transport dynamics, such as fluid flow in viscoelastic material [1], porous materials [2, 3], anomalous transport in biology [4], etc. Also fractional calculus is used to model problems in control theory, image processing, entropy theory and wave propagation phenomenon [5–7]. Due to large application in engineering, physics and mathematical sciences the analytical and numerical techniques for solving differential equations of fractional order increased dramatically in the recent year [8–11]. To learn more about fractional calculus, such as its existence, uniqueness of solution, applicability and solution methods, etc, [7, 12–14].

Many natural phenomena have been successfully modeled using fractional differential equations (FDEs) [7, 12, 15]. Physical phenomena including traffic flow models with fractional derivatives [16], fractional order KdV equation [9], fractional advection dispersion equation [10], and shock waves in a viscous medium [17] can be represented using non-linear differential equations.

The fractional derivative is a nonlocal operator, discretizing the temporal fractional derivative is a difficult task. To numerically solve FDEs in the complex domain and anticipate long-time ranges, an efficient numerical approach is required. Let's employ the Laplace transform to solve these problems because the temporal fractional derivatives are easily handled by the Laplace transform. Radial basis function generated finite difference (RBF-FD) method will be used to treat the spatial variable.

Mesh-free methods are becoming more popular, emerging, interesting and fascinating numerical techniques due to the ability to solve those physical and engineering problems with no meshing or minimum of meshing for which the traditionally used mesh-based methods are not suited like finite volumes, finite differences, finite elements, moving least square, point interpolation, element free Galerkin, reproducing kernel particle and boundary element free methods. RBFs methods appears to be really consists and most prominent techniques among the meshless methods family while looking at the interpolation of multi-dimensional scattered data and have received recently a tremendous and considerable attention in scientific community because of its capacity to achieve spectral accuracy, efficiency and high flexibility in solving complex PDEs, integral equations and fractional equations in comparison to other advanced approaches [18–21]. The most commonly used kernel in meshless techniques is the multi-quadric (MQ) kernel suggested by [21–24]. [22, 23] using the radial basis function for solution of a collocation strategy for PDEs.

Another specific qualitative characteristic disclosed on solutions to IBVPs of some evolutionary equations that have been established through investigations and are linked by their large-time stunts the term eventual time periodicity was coined to describe this phenomenon. This enticing and appealing event take place by a piston or flap or paddle-type wave maker put on one of the channel's ends in tests of research. When the wave generator oscillates at a predictable period $T_0 > 0$, at each location along the channel, the amplitude of the waves appears to become periodic, when a particular period of time has elapsed [25–28]. Various studies have previously addressed these important and fascinating eventual periodic phenomena such as Burger-type equations, generalized equations for KdV, BBM, and its dissipating counterparts [29–36].

In this work, let's offer an iterative strategy for obtaining approximate solutions of non-linear time fractional KdV and Burgers equation using the Laplace transform in conjunction with RBF-FD meshless method. Also the said scheme is implemented to investigate eventual periodicity to these equations.

## 2. Materials and methods
### 2. 1. Preliminaries
*Definition 1.*
The fractional order derivative in the Caputo sense [12, 37] is defined by:

$$D_t^\alpha w(t) = \frac{1}{\Gamma(n-\alpha)} \int_0^t \frac{1}{(t-z)^{\alpha+1-n}} \frac{d^n}{ds^n} w(s)\mathrm{d}s, \tag{1}$$

where

$$(n-1) < \alpha < n \in z^+.$$

*Definition 2.*
The Laplace transform of a given function $w(t)$ where $t \geq 0$ is defined as:

$$\hat{w}(t) = L\{w(t)\} = \int_0^t e^{-zt} w(t)\mathrm{d}t, \tag{2}$$

if improper integral converges.

*Lemma 1.*

Let $w(t)$ be a continuous function on $0 \le t \le t_n$ and if $M_1$, $M_2$ are some constants, with property,

$$\left| e^{-M_2 t} w(t) \right| < M_1, \ \forall t > t_n. \tag{3}$$

Then $Lw(t)$ exists.

*Lemma 2.*

If $w(t) \in C^p[0, \infty)$, then the Laplace transform of the Caputo fractional derivative is defined as:

$$L\left\{ D_t^\alpha w(t) \right\}(s) = s^\alpha \bar{w} - \sum_{i=0}^{n-1} s^{\alpha-i-1} s^{(i)}(0), \ \alpha \in (n-1, n) \in z^+. \tag{4}$$

## 2. 2. Nonlinear Time Fractional KdV-Burgers Model

A non-linear time fractional model for KdV-Burgers equation [38] is considered in the following form:

$$\frac{\partial^\alpha u(x,t)}{\partial t} + \eta u(x,t) \frac{\partial u(x,t)}{\partial x} - \xi \frac{\partial^2 u(x,t)}{\partial x^2} + \zeta \frac{\partial^3 u(x,t)}{\partial x^3} = f(x,t), \ x \in \Omega \subset R^d, \tag{5}$$

where $\alpha \in (0,1)$, $t \in [0,T]$, $\eta$, $\xi$, $\zeta$ are parameters and $f(x,t)$ is the source function.

$$u(x,0) = u_0, \ x \in \Omega \tag{6}$$

is the initial condition and

$$\mathcal{B} u(x,t) = g(t), \ x \in \partial\Omega, \ t \ge 0 \tag{7}$$

are the boundary conditions with boundary operator $\mathcal{B}$.

## 2. 3. Description of the method

*Description of RBF-Finite Differences method.*

Let's deal with general time dependent PDE for mathematical formulation and define the RBF-FD process in a gradual way. Take the problem of frame:

$$u_t(x,t) = \mathcal{L} u(x,t), \ x \in \Xi \subseteq R^d, \ d \ge 1, \ t > 0 \tag{8}$$

associated with initial and boundary conditions:

$$u(x,0) = u_0(x), \ \mathcal{B} u(x,t) = h(x,t), \ x \in \Xi, \tag{9}$$

where $u_0$ and $h$ are certain provided functions, while the spatial operators $\mathcal{L}$, $\mathcal{B}$ representing the differential operators. Assume $\{x_i\}_{i=1}^N$ denotes $N$ number of nodes used for approximation in the domain $\Xi$ for the given problem. RBF-FD is a mesh-free method and essentially a generalization of conventional finite difference (FD) method. In classical FD approach the derivative of a function $u$ is defined as a linear combination of the values of $u$ at some closest surrounding values (stencil) nodes. The difference is that RBF-FD methods use radial basis function instead of polynomials use in classical FD method [24].

*Global RBF differentiation matrix.*

Discretization of equations (8), (9) via global RBF method can be followed by approximating the unknown function $u$ by the linear combination of radial kernel $\phi$ at the node $x$ specified by:

$$\bar{u}(x) = \sum_{j=1}^N C_j \phi(\|x - x_j\|) = \Phi(x)^T C, \ x \in \Xi, \tag{10}$$

such that:

$$\Phi(x)^T = \left(\phi(\|x - x_1\|), \phi(\|x - x_2\|), ..., \phi(\|x - x_N\|)\right), \tag{11}$$

and $C$ is the expansion coefficients vector. Now equation (10) in Lagrange form is stated as:

$$\hat{u}(x) = \Phi(x)^T K^{-1} u. \tag{12}$$

Here $K$ representing system interpolation matrix for the global RBF. Now the interpolant (kernel-based) $\hat{u}$ in (12) gives good approximation of $u$. Consequently, any operator used on $\hat{u}$ also would be an excellent estimation of relevant operator employed on $u$ [21, 24]. Applying linear differential operator $\mathcal{L}$ on above equation (12) gives:

$$\mathcal{L}\hat{u}(x) = \mathcal{L}\Phi(x)^T K^{-1} u. \tag{13}$$

From (13) let's use the notation below for values:

$$K_{\mathcal{L}} = \begin{bmatrix} \mathcal{L}\Phi(x_1)^T \\ . \\ . \\ . \\ \mathcal{L}\Phi(x_N)^T \end{bmatrix}. \tag{14}$$

The global discretization (differentiation) matrix $L$ of size $N \times N$ may thus be considered as:

$$L = K_{\mathcal{L}} K^{-1}. \tag{15}$$

Since from (14), it is possible to see that the $i$th row of $K_{\mathcal{L}}$ corresponds to $\mathcal{L}\Phi(x_i)^T$, therefore let's observe from (15) that the $i$th row of $L$,

$$L_i = \mathcal{L}\Phi(x_i)^T K^{-1}, \tag{16}$$

serve as global differentiation matrix $L$ one single row.

*Local RBF differentiation matrix.*

Let's now report derivation of local differentiation matrix and describe how to compute the local finite differences associated weights which give rise to local interpolant in a locally small neighborhood regarding point $x_i$ exactly. Let's consider the set of points $\Xi = \{x_1, x_2, ..., x_N\}$, where it is possible to make the derivative approximation, these points can be regarded as stencil centers. For a given $i$th evaluation node say $x_i$, the size of nearest neighboring nodes in stencil $N_{x_i}$ of $x_i$ is $n$. Specifying also the set of points $Z = \{z_1, z_2, ..., z_N\}$ at which it is possible to analyze (sample) data. The points inside the stencil having size n are collected at $Z_i \subset Z$. Now estimation of differential operator $\mathcal{L}$ on stencil with center node $x_i$ and collected at $Z_i$ is given by:

$$L_i = K_{\mathcal{L}}^{x_i} K_{z_i}^{-1}. \tag{17}$$

Actually it assembles a stencil having center node $x_i$ hence to declare it as local differentiation matrix however it behaves globally since it operate whole entire data of that small stencil. All those $L_i$ matrices contains non-zeros entries in sparse (global) matrix $L^{FD}$, however their position must still be determined further in that sparse matrix $L^{FD}$.

Now $L_i^{FD}$ which representing the $i$th row of $L^{FD}$ and holds non-zero values from matrix $L_i$ (since it has one test node $x_i$ so it is row vector). As the points in $Z_i \subset Z$ are used in constructing $L_i$ (Hence columns of $L^{FD}$ connected along with those points which are non-zero columns of row $i$. Determining the position in the sparse row $L_i^{FD}$ of those points correctly, define an incidence matrix having entries below:

$$[P_i]_{k,l} = \begin{cases} 1, \text{if } k = l, i.e., k^{th} \text{ entry in } Z_i \text{ meet the } l^{th} \text{ in } Z, \\ 0, \text{ else.} \end{cases} \tag{18}$$

Use this to describe the complete sparse matrix as:

$$L^{FD} = \begin{bmatrix} K_L^{x_1} K_{z_1}^{-1} P_1 \\ . \\ . \\ . \\ K_L^{x_N} K_{z_N}^{-1} P_N \end{bmatrix}. \tag{19}$$

Ultimately the discretization for problem (8), (9) maybe written as:

$$u' = \mathcal{M} u, \tag{20}$$

where $\mathcal{M} = \begin{bmatrix} L^{FD} \\ B^{FD} \end{bmatrix}$, $B^{FD}$ stand for the discretization of operator applied at the boundary and can accordingly be found as $L^{FD}$. Evolving in time the ODE system (20), some solver ODE such as, ode113 ode23, ode45, and several others can be used from Matlab.

Time discretization via Laplace transform.

Now discuss how to reduce the non-linear fractional order differential equation (5)–(7) to system of nonlinear algebraic equations. The resulting system is then solved with the help of RBF-FD approach on a local level. The use of Laplace transform eliminates the time variable. Since the Laplace transform is a linear operator, it can't be directly applied to equations that are non-linear, therefore let's design an iterative approach for linearizing the nonlinear equation so that the Laplace transform can be applied to each step of the iterative process. It is possible to linearize the term $u(x,t)\dfrac{\partial u(x,t)}{\partial x}$ of the problem (5)–(7) as given below:

$$\frac{\partial^\alpha u^n(x,t)}{\partial t} + \eta u^{n-1}(x,t)\frac{\partial u^n(x,t)}{\partial x} - \xi \frac{\partial^2 u^n(x,t)}{\partial x^2} + \zeta \frac{\partial^3 u^n(x,t)}{\partial x^3} = f(x,t). \tag{21}$$

Applying Laplace transform to (21):

$$s^\alpha \bar{u}^n(x,t) - s^{\alpha-1}u_0 + \eta u^{n-1}(x,t)\frac{\partial \bar{u}^n(x,t)}{\partial x} - \xi \frac{\partial^2 \bar{u}^n(x,t)}{\partial x^2} + \zeta \frac{\partial^3 \bar{u}^n(x,t)}{\partial x^3} = \hat{f}(x,t), \tag{22}$$

and boundary conditions:

$$\mathcal{B}\,\bar{u}(x,s) = \hat{g}(s),\ x \in \partial\Omega. \tag{23}$$

Now introducing the differential operators $\mathcal{D}_1,\ \mathcal{D}_2$ and $\mathcal{D}_3$:

$$s^\alpha \bar{u}^n(x,s) - s^{\alpha-1}u_0 + \eta u^{n-1}(x,t)\mathcal{D}_1\bar{u}^n(x,s) - \xi \mathcal{D}_2\bar{u}^n(x,s) + \zeta \mathcal{D}_3\bar{u}^n(x,s) = \hat{f}(x,s). \tag{24}$$

Thus the non-linear iterative update process to approximate the solution $u(x,s)$, using $u_0$ as the initial approximation is given by:

$$\left[s^\alpha \mathcal{I} + \eta u^{n-1}(x,t)\mathcal{D}_1 - \xi \mathcal{D}_2 + \zeta \mathcal{D}_3\right]\bar{u}^n(x,s) = s^{\alpha-1}u_0 + \hat{f}(x,s),\ x \in \Omega \subset R^d,\ n = 1,2,..., \tag{25}$$

where $\bar{u}^n(x,s)$ is the $n^{\text{th}}$ iterate of the solution. The iterations will continue when the change in the solution is less than or equal to some pre-assign value. The problem (25) is approximated using radial basis functions in FD setting. An inverse Laplace transform approach is used for getting the solution $u(x,t)$ of the given model (21).

Technique of numerical inversion: to obtain the solution $u(x,t)$, let's use the numerical inversion formula below:

$$u(x,t) = \frac{1}{2\pi i} \int_{L_1}^{L_2} e^{st} \bar{u}(x,s) \mathrm{d}s, \ \lambda > \lambda_0, \tag{26}$$

where $L_1 = \lambda - i\infty$ and $L_2 = \lambda + i\infty$,

$$u(x,t) = \frac{1}{2\pi i} \oint e^{st} \bar{u}(x,s) \mathrm{d}s, \tag{27}$$

where $\Gamma$ is suitable path joining $L_1$ to $L_2$. In this paper let's use two types of contours. The parabolic contour [39]:

$$s = \nu(i\omega + 1)^2, \tag{28}$$

for the strip:

$$\omega = \varsigma + i\kappa, \tag{29}$$

where $-\infty < \varsigma < \infty$, which reduces to:

$$s(\varsigma) = \nu\left((1-\kappa)^2 - \varsigma^2\right) + 2i\nu\varsigma(1-\kappa), \tag{30}$$

and the Hyperbolic path [40] is defined as follows:

$$s(\varsigma) = \omega + \nu(1 - \sin(\alpha - i\varsigma)), \ \varsigma \in R, \tag{31}$$

with

$$\nu > 0, \ \omega \geq 0, \ 0 < \alpha < \lambda - \frac{1}{2}\pi$$

and

$$\frac{1}{2}\pi < \lambda < \pi \ [40].$$

Using the parabolic path defined in (30), integral in (26) becomes:

$$u(x,t) = \frac{1}{2\pi i} \int_{-\infty}^{\infty} e^{s(\varsigma)t} \bar{u}(x,s(\varsigma)) \bar{s}(\varsigma) \mathrm{d}\varsigma. \tag{32}$$

Applying the trapezoidal rule with equal step size the approximation to (32) is given by:

$$u_h(x,t) = \frac{h}{2\pi i} \sum_{j=-M}^{M} e^{s_j t} \bar{u}(x,s_j) s_j, \ s_j = s(\varsigma_j), \ \varsigma_j = jh, \tag{33}$$

where $h$ is the step size.

## 3. Results and Discussion

### 3. 1. Numerical Results of our proposed scheme

This section is devoted to demonstrate the accuracy and efficiency of our proposed numerical scheme called as radial basis function finite difference method combined with Laplace transform (RBF-FD)-LT by solving a number of nonlinear testing challenges. The main advantage of our method is avoidance of time stepping procedure which need a small time step for higher accuracy and stability, the second advantage is that our method has the capability to work in multi-dimensions in irregular domain as well. The accuracy of our method is tested in terms of $L_\infty$, $L_2$, root mean square (RMS) and $E(c, M)$ errors norms. Comparison with available results in literature for example in reference [38], **Tables 1–4** and **Fig. 1–9**, shows that our method is more accurate and stable. Let's solve some benchmark problems like time fractional KdV, Burger and KdV-Burger equations in one dimension as under.

Problem 1: consider model equation (5)–(7) for η = 1, ξ = 0, ζ = 1, let's obtain the non-linear KdV equation [41].

$$\frac{\partial^\alpha u(x,t)}{\partial t} + u(x,t)\frac{\partial u(x,t)}{\partial x} + \frac{\partial^3 u(x,t)}{\partial x^3} = f(x,t), \qquad (34)$$

where $0 < \alpha < 1$. The initial condition is:

$$u(x,0) = 0, \ x \in [a,b]$$

and

$$u(a,t) = \frac{1}{3000}e^{-a^2}t^5, \ t \geq 0; \ u(b,t) = \frac{1}{3000}e^{-b^2}t^5, \ t \geq 0,$$

are the boundary conditions:

$$f(x,t) = \frac{e^{-x^2}t^5}{25}\left(\frac{1}{\Gamma(6-\alpha)}t^{-\alpha} - \frac{1}{180000}t^2 x e^{-x^2} + \frac{1}{10}x - \frac{1}{15}x^3\right), \qquad (35)$$

the exact solution is:

$$u(x,t) = \frac{1}{3000}e^{-x^2}t^5. \qquad (36)$$

In this problem let's choose $[t_0, T] = [0.5, 5]$ and $t = 1$ in our computation.

The Hyperbolic contour (31) with optimal parameter values θ = 0.1, σ = 0.1541, ω = 0, is used and shown in **Fig. 1**.

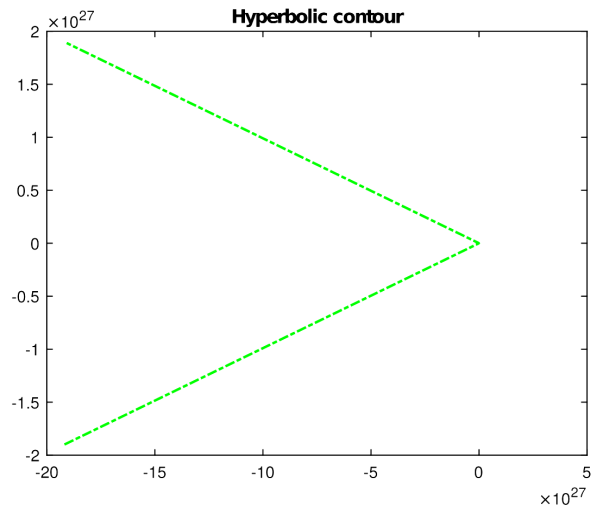

**Fig. 1.** Graphical representation of Hyperbolic contour (31) for $N = 71$, $M = 60$, $t_0 = 0.5$ and $T = 5$ corresponding to problem 1

The computed solutions at $N = 71$, $N_x = 5$, $M = 60$, $\alpha = 0.5$, $t = 1$ and $T = 5$ are shown in **Fig. 2**, our proposed numerical approach is stable and accurate, according to numerical computations.

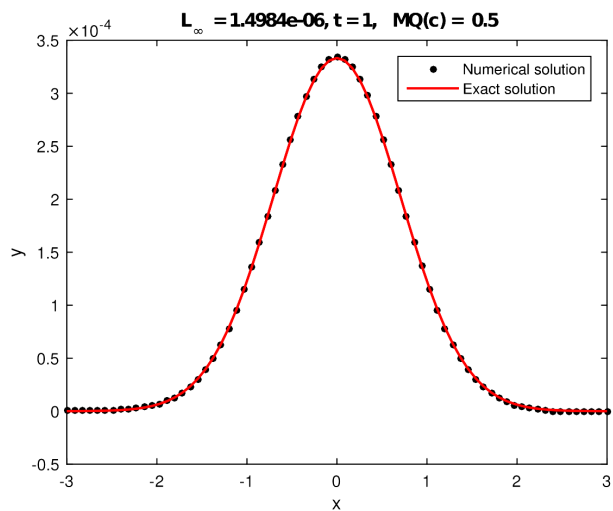

**Fig. 2.** Plot of approximate solution in comparison with exact solution (solid lines show exact solution and «.» showing numerical solution for problem 1

In **Table 1** the errors in terms of $L_\infty$, $L_2$ and RMS are decaying slowly for $21 \leq N \leq 71$, $M = 60$, $N_x = 5$, $\alpha = 0.5$, $[a, b] = [-3, 3]$, and the graphical representation of these errors are shown in **Fig. 3**. The estimation of error $E(c, M) = (e^{-cM/\log M})$ with $c = 0.5$ is $6.5742\text{E}^{-4}$.

**Table 1**
Numerical results corresponding to problem 1

| N | $L_\infty$ | $L_2$ | RMS |
|---|---|---|---|
| 21 | $1.0259\text{E}^{-05}$ | $2.6557\text{E}^{-05}$ | $4.8850\text{E}^{-07}$ |
| 41 | $3.0939\text{E}^{-06}$ | $1.1907\text{E}^{-05}$ | $7.5460\text{E}^{-08}$ |
| 51 | $3.5167\text{E}^{-06}$ | $1.5228\text{E}^{-05}$ | $6.8955\text{E}^{-08}$ |
| 61 | $1.0960\text{E}^{-06}$ | $4.4012\text{E}^{-06}$ | $1.7968\text{E}^{-08}$ |
| 71 | $1.4984\text{E}^{-06}$ | $6.7997\text{E}^{-06}$ | $2.1104\text{E}^{-08}$ |
| 81 | $6.7697\text{E}^{-07}$ | $2.9525\text{E}^{-06}$ | $8.3576\text{E}^{-09}$ |
| 91 | $9.0620\text{E}^{-07}$ | $4.7338\text{E}^{-06}$ | $9.9583\text{E}^{-09}$ |
| 100 | $7.6117\text{E}^{-07}$ | $4.3381\text{E}^{-06}$ | $7.6117\text{E}^{-09}$ |

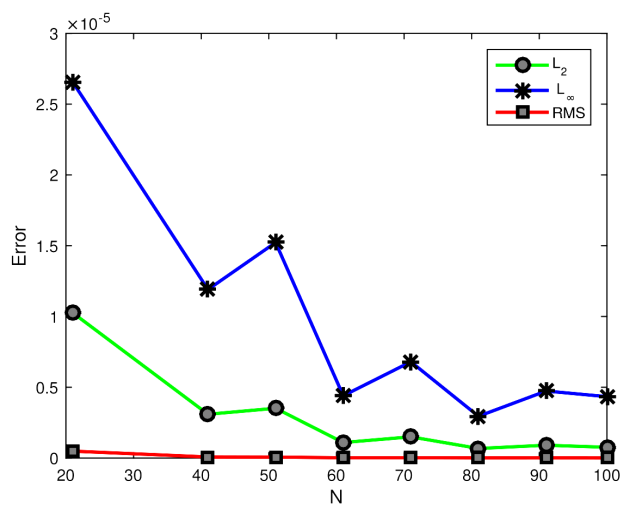

**Fig. 3.** Graphical representation of errors obtained in **Table 1** corresponding to problem 1

**Table 2** shows the $L_\infty$ error and error estimate $E(c,M)$ for various α values and quadrature points along the hyperbolic contour (31).

**Table 2**
Errors corresponding to problem 1

| $(N, N_x = 71.5)$ | $\alpha = 0.20$ | $\alpha = 0.50$ | $\alpha = 0.75$ | $\alpha = 0.90$ | $E(c, M)$ |
|---|---|---|---|---|---|
| $M$ | $L_\infty$ | $L_\infty$ | $L_\infty$ | $L_\infty$ | |
| 45 | 1.6109E$^{-06}$ | 1.5257E$^{-06}$ | 1.4358E$^{-06}$ | 1.4302E$^{-06}$ | 2.7103E$^{-03}$ |
| 61 | 1.6014E$^{-06}$ | 1.4984E$^{-06}$ | 1.3521E$^{-06}$ | 1.2451E$^{-06}$ | 5.9954E$^{-04}$ |
| 71 | 1.6014E$^{-06}$ | 1.4984E$^{-06}$ | 1.3521E$^{-06}$ | 1.2451E$^{-06}$ | 2.4163E$^{-04}$ |
| 91 | 1.6014E$^{-06}$ | 1.4984E$^{-06}$ | 1.3521E$^{-06}$ | 1.2451E$^{-06}$ | 4.1627E$^{-05}$ |

Problem 2: now consider model equation (5)–(7) for η = 1, ξ = 1, ζ = 0, let's obtain the non-linear Burger equation [38] with:

$$f(x,t) = \frac{2t^{2-\alpha}\sin(2\pi x)}{\Gamma(3-\alpha)} + 2\pi t^4 \sin(2\pi x)\cos(2\pi x) + 4\mu\pi^2 t^2\sin(2\pi x), \qquad (37)$$

$$u(x,0) = 0, \qquad (38)$$

$$u(0,t) = 0, \; u(1,t) = 0, \; t \geq 0, \qquad (39)$$

the actual solution is:

$$u(x,t) = t^2\sin(2\pi x). \qquad (40)$$

Here let's consider the domain $\Omega = [0, 1]$. The contour mentioned in (30) is used here and is shown in **Fig. 4**.

Various numerical solutions are shown in **Tables 3, 4**. The computed solutions are almost converges to exact solutions, as seen in **Table 3**, where let's choose $[t_0,T] = [0.09995,3]$, $t = 0.1$, $\alpha = 0.5$, $N = 81 \in \Omega$, stencils points $N_x = 9 \in \Omega$, $\xi = 0.10, 0.5, 1$ and quadrature points $M = 80, 90, 100$. At the same values of parameters ξ, $t$, α and ζ as given in[38] the proposed techniques yields identical results to those reported in [38]. Finally, the numerical solution at $\xi = 1$, $\alpha = 0.5$, $N = 41$, $N_x = 5$ stencil points $t_0 = 0.3$, $t = 0.5$, $T = 3$ and $M = 50$ quadrature points are depicted in **Fig. 5**. The proposed numerical scheme is accurate and stable.

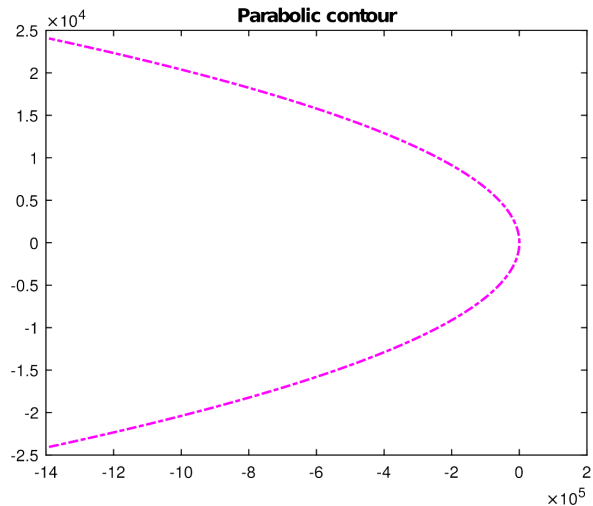

**Fig. 4.** Graphical representation of parabolic contour (30) for $N = 81$, $M = 60$, $t_0 = 0.09995$, $t = 0.1$ and $T = 3$ corresponding to problem 2

**Table 3**

Comparison table for approximate vs actual solution corresponding to problem 2

| $x$ | $M = 80$ | $M = 90$ | $M = 100$ | Exact [38] |
|---|---|---|---|---|
| | $\xi = 1$ | $\xi = 0.5$ | $\xi = 0.1$ | |
| 0.0 | 0.0000 | 0.0000 | 0.0000 | 0.000000 |
| 0.1 | 0.0059 | 0.0059 | 0.0059 | 0.005878 |
| 0.2 | 0.0095 | 0.0095 | 0.0095 | 0.009511 |
| 0,3 | 0.0095 | 0.0095 | 0.0095 | 0.009511 |
| 0.4 | 0.0059 | 0.0059 | 0.0059 | 0.005878 |
| 0.5 | 0.0000 | 0.0000 | 0.0000 | 0.000000 |
| 0.6 | −0.0059 | −0.0059 | −0.0059 | −0.005878 |
| 0.7 | −0.0095 | −0.0095 | −0.0095 | −0.009511 |
| 0.8 | −0.0095 | −0.0095 | −0.0095 | −0.009511 |
| 0.9 | −0.0059 | −0.0059 | −0.0059 | −0.005878 |
| 1 | 0.0000 | 0.0000 | 0.0000 | 0.000000 |
| (RBF-FD)-LT | $L_\infty = 3.9383\text{E}^{-6}$ | $L_\infty = 4.0820\text{E}^{-6}$ | $L_\infty = 1.4288\text{E}^{-5}$ | |
| [38] | $L_\infty = 4.1294\text{E}^{-5}$ | $L_\infty = 3.7739\text{E}^{-5}$ | $L_\infty = 2.1269\text{E}^{-5}$ | – |
| (RBF-FD)-LT | $L_2 = 2.5970\text{E}^{-5}$ | $L_2 = 2.7051\text{E}^{-5}$ | $L_2 = 1.0798\text{E}^{-4}$ | |
| [38] | $L_2 = 2.9174\text{E}^{-5}$ | $L_2 = 2.6666\text{E}^{-3}$ | $L_2 = 1.5070\text{E}^{-5}$ | |

**Table 4**

Errors at $\xi = 1$, $t_0 = 0.3$, $t = 0.5$ and $T = 5$ corresponding to problem 2

| $(N, N_x = 61.11)$ | $\alpha = 0.2$ | $\alpha = 0.50$ | $\alpha = 0.75$ | $\alpha = 0.90$ | $\alpha = 0.5$ |
|---|---|---|---|---|---|
| $M$ | $L_\infty$ | $L_\infty$ | $L_\infty$ | $L_\infty$ | RMS |
| 40 | $1.2397\text{E}^{-03}$ | $1.2292\text{E}^{-03}$ | $1.2145\text{E}^{-03}$ | $1.2018\text{E}^{-03}$ | $2.0151\text{E}^{-05}$ |
| 50 | $1.2278\text{E}^{-03}$ | $1.2172\text{E}^{-03}$ | $1.2025\text{E}^{-03}$ | $1.1898\text{E}^{-03}$ | $1.9954\text{E}^{-05}$ |
| 60 | $1.2277\text{E}^{-03}$ | $1.2172\text{E}^{-03}$ | $1.2024\text{E}^{-03}$ | $1.1897\text{E}^{-03}$ | $1.9954\text{E}^{-05}$ |
| 70 | $1.2277\text{E}^{-03}$ | $1.2172\text{E}^{-03}$ | $1.2024\text{E}^{-03}$ | $1.1897\text{E}^{-03}$ | $1.9954\text{E}^{-05}$ |
| 80 | $1.2277\text{E}^{-03}$ | $1.2172\text{E}^{-03}$ | $1.2024\text{E}^{-03}$ | $1.1897\text{E}^{-03}$ | $1.9954\text{E}^{-05}$ |
| $(M, N_x = 50.5)$ | $\alpha = 0.2$ | $\alpha = 0.50$ | $\alpha = 0.75$ | $\alpha = 0.90$ | $\alpha = 0.5$ |
| $N$ | $L_\infty$ | $L_\infty$ | $L_\infty$ | $L_\infty$ | RMS |
| 11 | $1.8568\text{E}^{-03}$ | $1.8389\text{E}^{-03}$ | $1.8155\text{E}^{-03}$ | $1.7963\text{E}^{-03}$ | $1.6717\text{E}^{-04}$ |
| 21 | $1.1852\text{E}^{-03}$ | $1.1750\text{E}^{-03}$ | $1.1608\text{E}^{-03}$ | $1.1486\text{E}^{-03}$ | $5.5954\text{E}^{-05}$ |
| 41 | $1.2308\text{E}^{-03}$ | $1.2203\text{E}^{-03}$ | $1.2055\text{E}^{-03}$ | $1.1928\text{E}^{-03}$ | $2.9762\text{E}^{-05}$ |
| 51 | $1.2304\text{E}^{-03}$ | $1.2198\text{E}^{-03}$ | $1.2050\text{E}^{-03}$ | $1.1923\text{E}^{-03}$ | $2.3918\text{E}^{-05}$ |
| 61 | $1.2278\text{E}^{-03}$ | $1.2172\text{E}^{-03}$ | $1.2025\text{E}^{-03}$ | $1.1898\text{E}^{-03}$ | $1.9954\text{E}^{-05}$ |

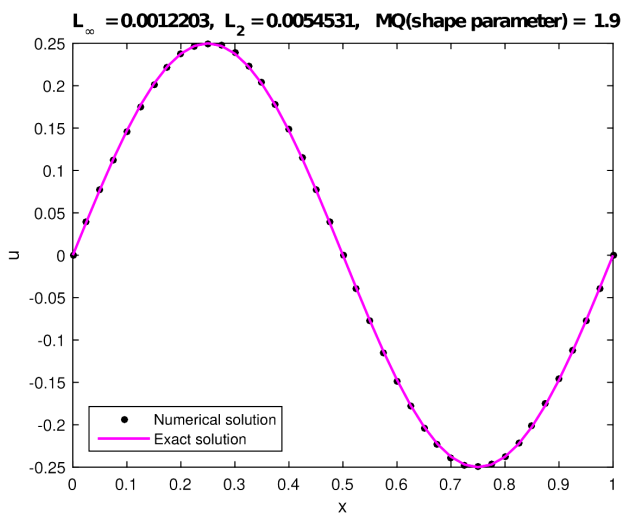

**Fig. 5.** Plot of approximate solution for problem 2 at $N = 41$, $N_x = 5$, $\alpha = 0.5$, $\xi = 1$, $M = 50$, $t_0 = 0.3$, $t = 0.5$ and $T = 3$ in comparison with exact solution (solid lines show exact solution and «.» showing numerical solution

Problem 3: finally, look at the basic Burger equation [42] which is given by:

$$\frac{\partial^{\alpha}u(x,t)}{\partial t} + \eta u(x,t)\frac{\partial u(x,t)}{\partial x} - \xi\frac{\partial^2 u(x,t)}{\partial x^2} = 0. \tag{41}$$

The initial condition is:

$$u(x,0) = \frac{2\pi\xi\sin(\pi x)}{\eta + \cos(\pi x)}, \ \eta > 1, \ \xi > 0, \tag{42}$$

and

$$u(0,t) = 0, \ u(1,t) = 0, \ t > 0 \tag{43}$$

are the boundary conditions.

The exact solution is:

$$u(x,t) = \frac{2\pi\xi\exp(-\pi^2 t)\sin(\pi x)}{\eta + \exp(-\pi^2\xi t)\cos(\pi x)}, \ \eta > 1. \tag{44}$$

Here take $\Omega = [0,1]$. The parabolic contour (30) is used here. The parameter values $t = 0.001$, $\xi = 0.1, 0.2, 0.5, 1$, and $N = 41$, and as given in [42] are used. **Fig. 6** depicts the numerical solutions found using our proposed method, the results of the numerical technique are identical to those in [42].

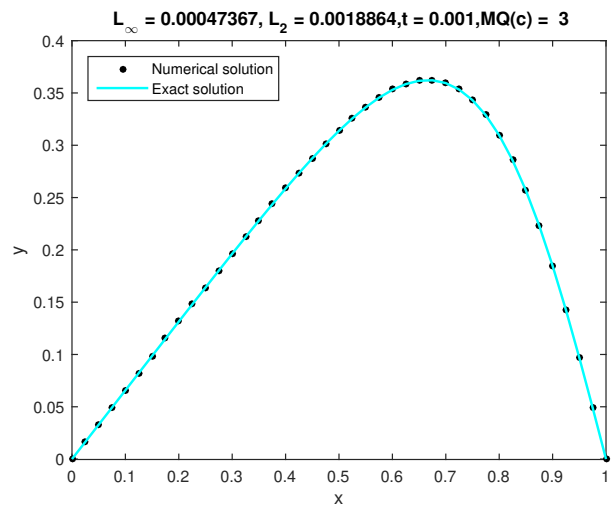

**Fig. 6.** Numerical result for $N = 41$, $M = 71$, $t_0 = 0.0001$, $T = 1$, $\eta = 2$ and $\xi = 0.1$ corresponding to problem 3

### 3. 2. Eventual Periodicity of the proposed model equations

Now let's present the results of our method investigating the eventual periodicity of the given model equation (5)–(7) for KdV, Burger and KdV-Burger equations studied by [30] in graphical form along with appropriate boundary data $g(t)$ (periodic function of period $T_0 > 0$). The initial data $u_0$ is not necessary in eventual periodicity so take it zero. For each problem the amplitudes $u(x,t)$ produced in six graphs at particular points. $N$ indicates complete domain points, while $N_x$ denotes points in respective sub-domains. The $X$ and $Y$ axes are representative in these graphs of time $t$ and amplitude $u$ respectively. The last graph shows the amplitude remains zero in every problem.

Eventual periodicity of fractional KdV equation: let's compute eventual periodicity of model equation (5)–(7) for fractional KdV equation with parameters $\alpha = 0.2$, $\eta = 0.05$, $\xi = 0$, $\zeta = 1$, $f(x,t) = 0$ and $g(t) = \sin(20\pi t)\tanh(5t)$. The amplitudes $u(x,t)$ for this model is shown in six plots in **Fig. 7** at given specific points $x = -0.950670$, $-0.808460$, $-0.587280$, $-0.308720$, $0.0$, $0.999650$ in the domain $[-1, 1]$ and in a time domain $[0, 1.8]$. The plots below clearly confirm the subsequent periodic activity with damped amplitude of the solution in the specified domain at these particular positions.

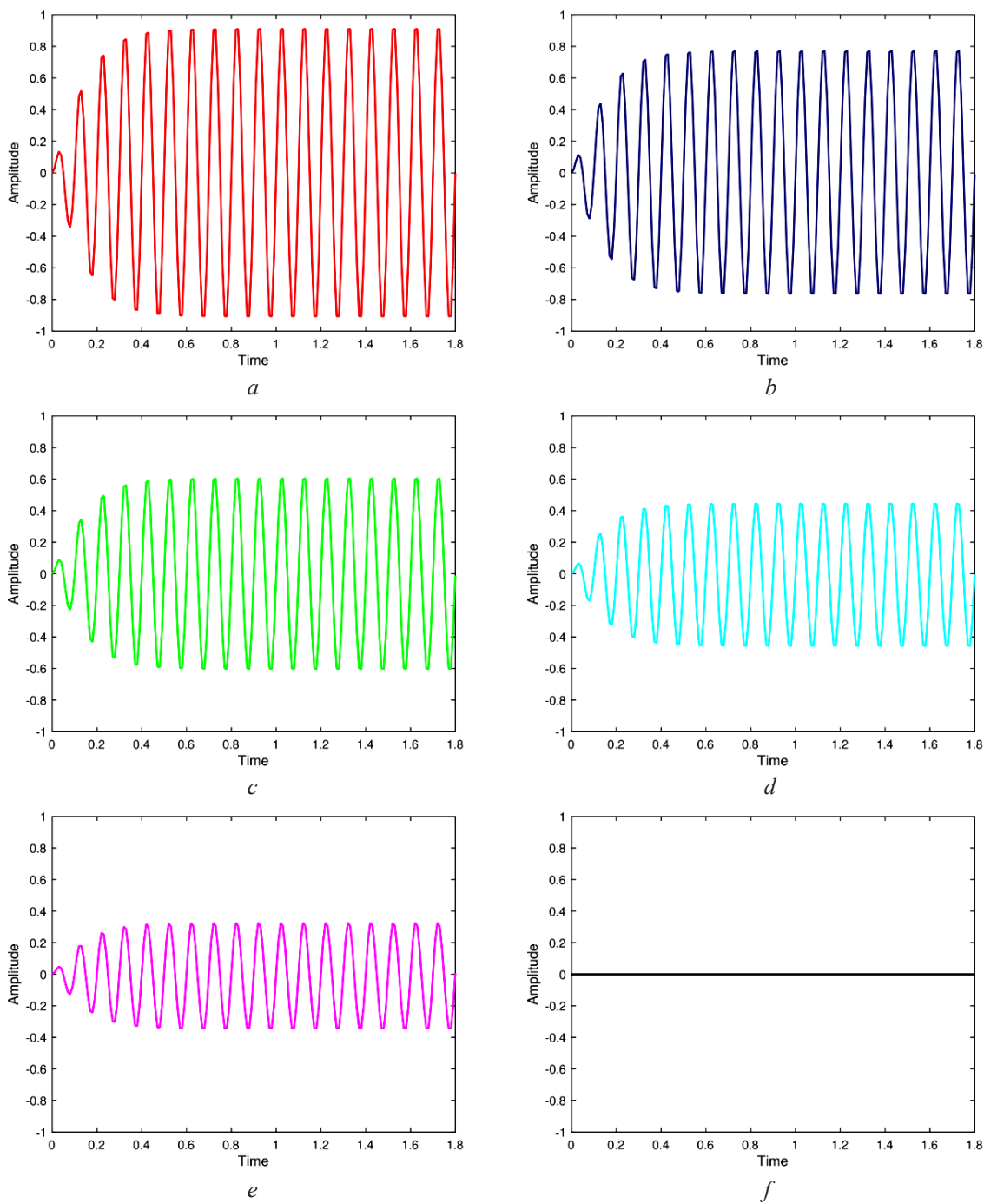

**Fig. 7.** Eventual periodicity for fractional KdV equation at different specific points in space domain [−1, 1] and in a time domain [0, 1.8] as shown in above plots: $a$ − at $x = −0.950670$; $b$ − at $x = −0.808460$; $c$ − at $x = −0.587280$; $d$ − at $x = −0.308720$; $e$ − at $x = 0.0$; $f$ − at $x = 0.999650$, using $N = 200$, $N_x = 25$, $g(t) = \sin(20\pi t)\tanh(5t)$

Eventual periodicity of fractional Burger equation: compute eventual periodicity of model equation (5)–(7) for fractional Burger equation with parameters $\alpha = 0.2$, $\eta = 0.05$, $\xi = 1$, $\zeta = 0$, $f(x,t) = 0$ and $g(t) = \sin(20\pi t)\tanh(5t)$. The amplitudes $u(x,t)$ for this model is shown in six plots

144

in **Fig. 8** at given specific points $x = -0.950670$, $-0.808460$, $-0.587280$, $-0.308720$, $0.0$, $0.999650$ in the domain $[-1, 1]$ and in a time domain $[0, 1.8]$. The plots below clearly confirm the subsequent periodic activity of the solution in the specified domain at these particular positions. The influence of Burgers type dissipation on the damped amplitude is clearly seen.

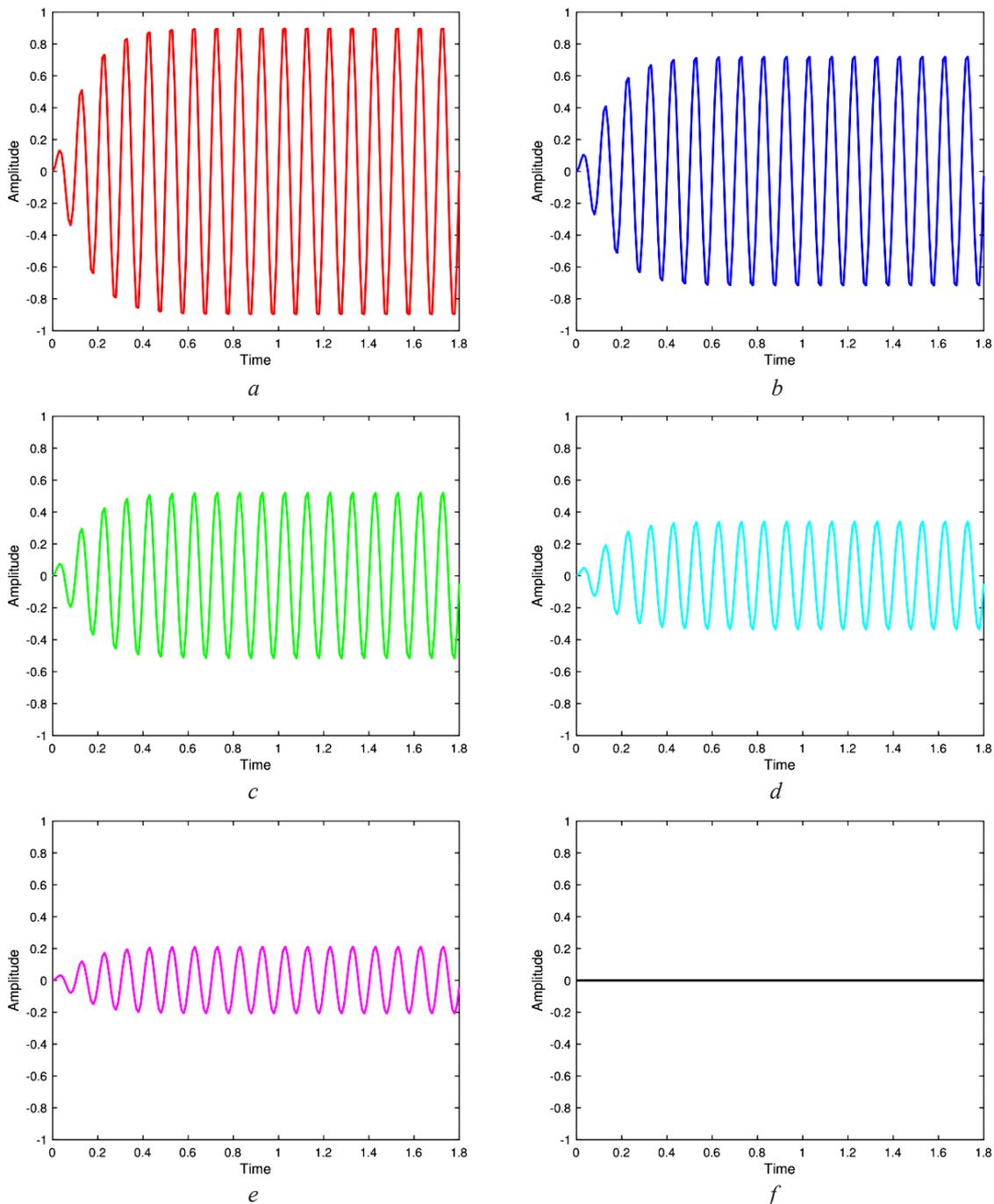

**Fig. 8.** Eventual periodicity for fractional Burger equation at different specific points in space domain $[-1, 1]$ and in a time domain $[0, 1.8]$ as shown in above plots $a$ − at $x = -0.950670$; $b$ − at $x = -0.808460$; $c$ − at $x = -0.587280$; $d$ − at $x = -0.308720$; $e$ − at $x = 0.0$; $f$ − at $x = 0.999650$, using $N = 200$, $N_x = 25$, $g(t) = \sin(20\pi t)\tan h(5t)$

Eventual periodicity of fraction KdV-Burger equation: finally compute eventual periodicity of model (5)–(7) for fractional KdV-Burger equation with parameters $\alpha = 0.2$, $\eta = 0.05$, $\xi = 10^{-4}$, $\zeta = 10^{-5}$, $f(x,t) = 0$ and $g(t) = \sin(20\pi t)\tanh(5t)$. The amplitudes $u(x,t)$ for this model is shown in six plots in **Fig. 9** at given particular points $x = -0.950670$, $-0.808460$, $-0.587280$, $-0.308720$, $0.0$, $0.999650$ in the domain $[-1, 1]$ and in a time domain $[0, 1.8]$. The plots below clearly confirm the subsequent periodic activity of the solution in the specified domain at these particular positions. The pattern of eventual periodicity has not changed, but the amplitudes have been drastically reduced.

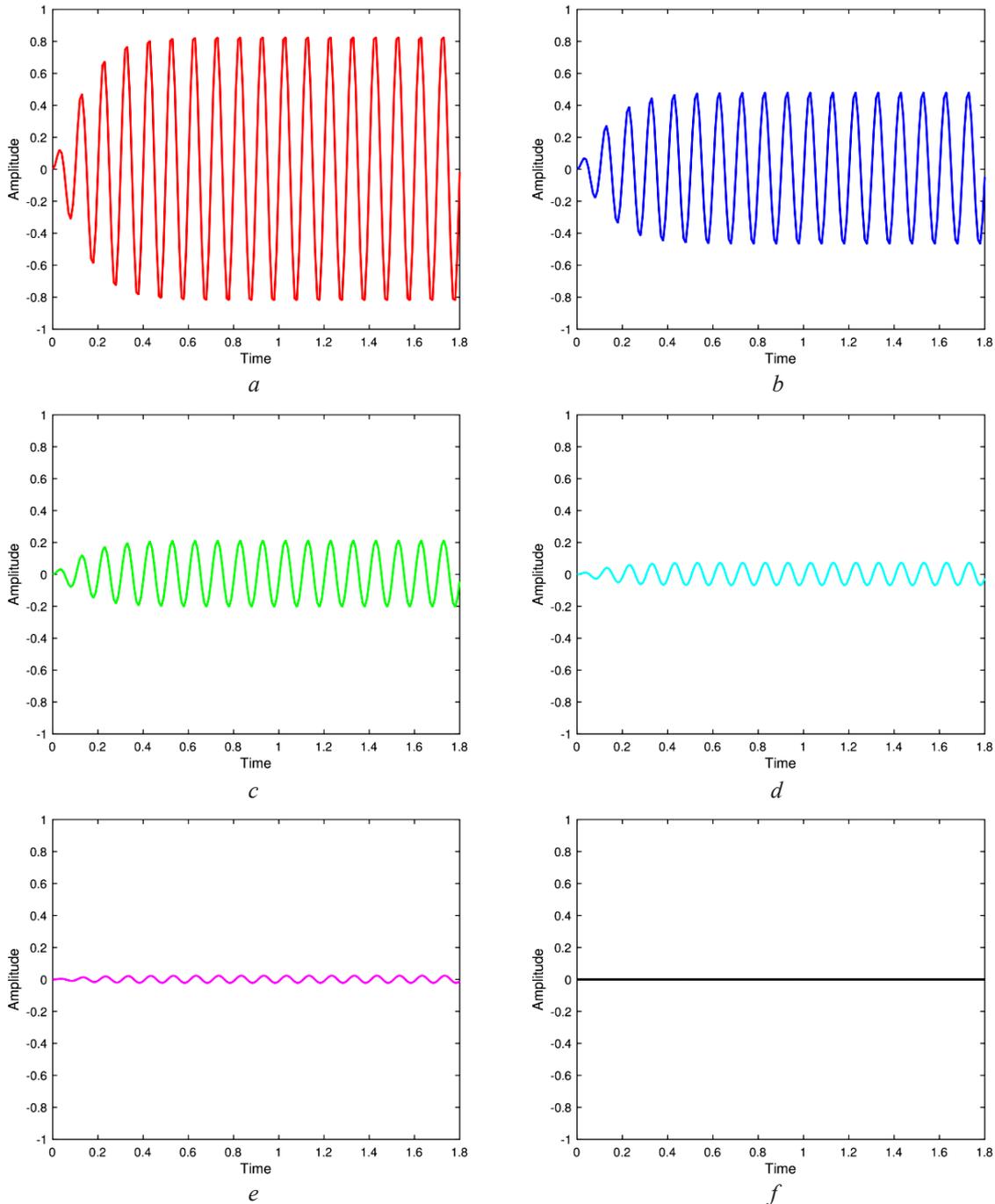

**Fig. 9.** Eventual periodicity for fractional KdV-Burger equation at different specific points in space domain $[-1, 1]$ and in a time domain $[0, 1.8]$ as shown in above plots:
$a$ – at $x = -0.950670$; $b$ – at $x = -0.808460$; $c$ – at $x = -0.587280$;
$d$ – at $x = -0.308720$; $e$ – at $x = 0.0$; $f$ – at $x = 0.999650$,
using $N = 200$, $N_x = 25$, $g(t) = \sin(20\pi t)\tanh(5t)$

146

## 4. Conclusions

In this research work, a hybrid numerical scheme is developed by combining the Laplace transform technique with RBF-FD mesh free method for approximating the solution of some non-linear time fractional dispersive wave equations (like KdV, Burger and KdV-Burger model) and also the periodic behavior of solution called as eventual periodicity. This method is integrated with Laplace transform approach for time integration and for fractional derivative in Caputo sense. The main advantage of Laplace transform technique is that it is free of issues relating to stability due to time steeping scheme. The spatial operators in multi-dimensions are approximated by RBF in the finite difference (FD) setting which generates small size differentiation matrices in local sub-domains and these are assembled as a single sparse matrix in the global domain. So large amount of data can be manipulated very easily and accurately. The construction of our approach is simpler and easier to solve any nonlinear higher order time fractional PDEs as compared to other numerical methods available in the literature. The simulation has been demonstrated through error norms recorded in tables and graphs, which shows the accuracy and robustness of our method.